\documentclass[10pt]{amsart}

\usepackage{amsmath, amssymb, amsfonts,
amsthm, amscd, latexsym, enumerate}
\theoremstyle{plain}
\newtheorem{theorem}{Theorem}[subsection]
\newtheorem{lemma}[theorem]{Lemma}
\newtheorem{proposition}[theorem]{Proposition}
\newtheorem{corollary}[theorem]{Corollary}

\theoremstyle{definition}
\newtheorem{definition}[theorem]{Definition}
\newtheorem{example}[theorem]{Example}
\newtheorem{remark}[theorem]{Remark}

\numberwithin{equation}{subsection}

\newcommand\fr{\mathfrak}
\renewcommand\sc{\mathcal}
\newcommand\scr{\mathcal}

\renewcommand{\r}[1]{~\ref{#1}}
\newcommand\bref{\eqref}

\newcommand{\CBbb}{\mathbb C}
\newcommand{\EBbb}{\mathbb E}

\newcommand{\WBbb}{\mathbb W}

\newcommand{\A}{\mathcal A}
\newcommand{\Ce}{\mathcal C}
\newcommand{\E}{\mathcal E}
\newcommand{\G}{\mathcal G}

\renewcommand{\O}{\Omega}
\renewcommand{\o}{\omega}

\newcommand\dimc{\dim_{\CBbb}}

\renewcommand{\a}{\alpha}
\newcommand{\be}{\beta}

\newcommand{\Ga}{\Gamma}

\newcommand{\si}{\sigma}

\newcommand{\ra}{\rightarrow}

\newcommand{\lra}{\longrightarrow}

\newcommand\hra{\hookrightarrow}

\newcommand\bra{\mapsto}

\newcommand{\oset}[1]{\overset {#1}{\ra}}
\newcommand{\osetl}[1]{\overset {#1}{\lra}}

\newcommand{\ocong}[1]{\overset {#1}  {\cong}}

\newcommand\ot{\otimes}

\newcommand{\tp}[2]{{#1} \ot {#2}}

\newcommand{\tpc}[2]{{#1} \ot_{\mathbb C} {#2}}

\newcommand{\exts}[5]{\Ext^1_{{\scr O}_{#5}}({\scr #1}_{#2} , {\scr #3}_{#4})}

\newcommand{\schom}[2]{{\sc H}om_{\mathbb C}( {#1} , {#2} )}

\newcommand{\scend}[1]{{\sc E}nd_{\mathbb C}( {#1} )}

\newcommand\del{\partial}
\newcommand\delbar{\bar \del}

\newcommand\ti{\tilde}
\newcommand\wti{\widetilde}

\newcommand\vphi{\varphi}

\newcommand\Id{\operatorname{Id}}

\newcommand\Aut{\operatorname{Aut}}

\newcommand\Hom{\operatorname{Hom}}
\newcommand\End{\operatorname{End}}

\newcommand\im{\operatorname{im}}

\newcommand\Lie{\operatorname{Lie}}

\newcommand\gl{\operatorname{\mathfrak {gl}}}

\newcommand\Tr{\operatorname{Tr}}
\newcommand\Ext{\operatorname{Ext}}
\newcommand\Iso{\operatorname{Iso}}
\newcommand\Hol{\operatorname{Hol}}
\newcommand\Holiso{\operatorname{Hol}_{\Iso}}
\newcommand\Vol{\operatorname{Vol}}

\newcommand\dvol{\operatorname{dvol}}

\newcommand\rank{\operatorname{rank}}

\newcommand\Dirac{{\mathpalette\c@ncel D}}
\newcommand\dirac{{\mathpalette\c@ncel \del}}
\newcommand\caldir{{\mathpalette\c@ncel \De}}

\newcommand{\extn}{0 \ra {\sc E}_1 \lra {\sc E} \lra {\sc E}_2 \ra 0}

\newcommand{\emc}[1]{\End _{\mathbb C}({#1})}

\newcommand{\emg}[2]{\End_{#2}( {#1} )}
\newcommand{\homg}[3]{\Hom_{#3}( {#1} , {#2} )}

\newcommand{\summa}[3]{\sum_{#1}^{#2} ~{#3}}

\newcommand{\ov}[2]{ \frac{#1}{#2} }

\newcommand{\degx}[2]{\deg_{#2} (#1)}
\newcommand{\degs}[1]{\degx{#1}{\sigma}}

\newcommand{\slox}[2]{\mu_{#1} (#2)}
\newcommand{\slos}[1]{\slox{#1}{\si}}

\newcommand{\lam}[1]{\mu_{#1}}

\newcommand{\lambdax}[1]{\Lambda_{#1}}
\newcommand{\lambdas}{\lambdax{\sigma}}

\newcommand{\lamj}[2]{\lam{{#1}_{#2}}}
\newcommand{\sclamj}[2]{\lamj{\sc #1}{#2}}

\newcommand{\fibint}{\int_F ~\Tr}

\newcommand{\ldyncbbb}[4]{\begin{picture}(75,12)
                       \put(0,2){\line(5,0){20}}
                       \put(20,2){\line(5,0){20}}
                       \put(45,2){\line(5,0){20}}
                       \put(-4,-1){$\times$}
                       \put(20,2){\circle*{4}}
                       \put(45,2){\circle*{4}}
                       \put(65,2){\circle*{4}}
                       \put(-2,6){\shortstack{\small #1}}
                       \put(17,6){\shortstack{\small #2}}
                       \put(43,6){\shortstack{\small #3}}
                       \put(62,6){\shortstack{\small #4}}
                       \end{picture}}

\newcommand{\ldyncbcc}[4]{\begin{picture}(75,12)
                       \put(0,2){\line(5,0){20}}
                       \put(20,2){\line(5,0){20}}
                       \put(45,2){\line(5,0){20}}
                       \put(-4,-1){$\times$}
                       \put(20,2){\circle*{4}}
                       \put(41,-1){$\times$}
                       \put(61,-1){$\times$}
                       \put(-2,6){\shortstack{\small #1}}
                       \put(17,6){\shortstack{\small #2}}
                       \put(43,6){\shortstack{\small #3}}
                       \put(62,6){\shortstack{\small #4}}
                       \end{picture}}

\begin{document}

\noindent
Proc. in Global Anaysis, Diff. Geometry and Lie Algebras, 
\newline\noindent
Balkan Geometry Press, Bucharest 2001, pp. 43--62. 

\vskip 1in

\title[Dimensional reduction of the perturbed Hermitian--Einstein equation]
{Dimensional reduction of   \\
the perturbed Hermitian--Einstein equation}

\date{\today}

\author[S. B.~Bradlow]{Steven B.~Bradlow${}^1$}

\address{Department of Mathematics\\
University of Illinois at Urbana--Champaign \\
Urbana IL 61801 USA}

\email[S. B.~Bradlow]{bradlow@math.uiuc.edu}

\thanks{${}^1$ Supported in part by The National Science
Foundation under Grant DMS-9703869.}

\author[J. F.~Glazebrook]{James F.~Glazebrook}

\address{Department of Mathematics\\
Eastern Illinois University\\
Charleston IL 61920 USA \\
and Department of Mathematics \\
University of Illinois at Urbana--Champaign \\
Urbana IL 61801 USA}

\email[J. F.~Glazebrook]{cfjfg@eiu.edu , glazebro@math.uiuc.edu}

\author[F. W.~Kamber]{Franz W.~Kamber${}^2$}

\address{Department of Mathematics\\
University of Illinois at Urbana--Champaign \\
Urbana IL 61801 USA}

\email[F. W.~Kamber]{kamber@math.uiuc.edu}

\thanks{${}^2$Supported in part by The National Science Foundation under
Grant DMS-9504084.}


\begin{abstract}
Given a K\"ahlerian holomorphic fiber bundle $F \hra M \ra X$,
whose fiber $F$ is a compact homogeneous K\"ahler manifold, we
describe the perturbed Hermitian--Einstein equations relative to
certain holomorphic vector bundles ${\sc E} \ra M$~. With respect
to special metrics on $\sc E$, there is a dimensional reduction
procedure which reduces this equation to a system of equations on
$X$ known as the twisted coupled vortex equations.
\end{abstract}

\maketitle



\section{Introduction}

Dimensional reduction techniques are applicable to studying
special solutions to partial differential equations particularly
in the presence of a group action where invariant solutions are of
interest. The invariant solutions may be interpreted as solutions
to an associated set of equations on a lower dimensional space of
orbits of the group action. However, one may ask if there is no
group action, is it still possible to dimensionally--reduce the
original system? A positive answer points to the study of the
Hermitian--Einstein (HE) equation with respect to special metrics
on holomorphic bundles together with some extra data. The purpose
of this paper is to outline a construction leading to dimensional
reduction of a class of equations which we call the {\it perturbed
Hermitian--Einstein equations} (briefly, the PHE equations) on a
Hermitian holomorphic vector bundle ${\sc E} \lra M$ where $M$ is a
compact K\"ahler manifold. We stress that the term perturbed
has here a delicate interpretation as will be apparent from the
text. In fact, the PHE equations are actually more general than
the HE equations because they possess an extra perturbation term.
This extra term arises from the fact that in this case,
$M$ is the total space of a holomorphic fiber bundle $F \hra
M \lra X$, where  $X$ is a compact K\"ahler manifold and the fiber
$F$ is a compact K\"ahlerian homogeneous space. Now $\E$ as a holomorphic
vector bundle is obtained via a holomorphic extension of certain
holomorphic vector bundles on $M$ and is equipped with an
invariant hermitian metric. This metric together with
the Kobayashi form of the extension and some natural conditions
on $F$, imply that the PHE equation is equivalent to a system of
equations on $X$, namely the {\it twisted coupled vortex
equations}.

\medbreak
The overall construction, on which there are several variations,
relies on results relating to the representation theory
of complex semisimple Lie groups and the
Bott--Borel--Weil theorem. In addition, the PHE equation can be
obtained as a moment map equation. Here we will outline the
general construction of \cite{BGKfour} leading to the twisted coupled
vortex equations (cf \cite{BGKone} \cite{BGKtwo} \cite{GProne}).
The existence theory of the solutions of such twisted coupled
equations is discussed
via the Hitchin--Kobayashi correspondence in
\cite{BGKthree}. References \cite{AGone} \cite{AGtwo} contain an
independent study of several aspects of this theory and focus on
other questions.


\section{Some preliminaries}

\subsection{The K\"ahler manifold $M$}

Let us commence by describing the compact homogeneous K\"ahler manifold $F$,
that is, for connected complex Lie groups $G$ and $P$ with $G$ semisimple
and $P \subset G$ parabolic, we set $F = G / P \cong U / K$ where
\begin{equation}\label{symmetric}
G = \Hol (F)_e~, ~U = \Holiso (F)_e~, ~K = U \cap P~.
\end{equation}
Furthermore, $F$ is a simply connected algebraic manifold, the groups
$U$ and $K$ are connected compact Lie groups, with $U$ semisimple and
$K$ the centralizer of a torus (hence $K \subset U$ has maximal
rank) and any $G$--invariant hermitian metric on $F$ is a K\"ahler
(for further details see \cite{BE} \cite{Botttwo} \cite{Kobtwo}).
The equivariant holomorphic vector bundles on $G / P$ are
homogeneous vector bundles \cite{Botttwo} given by
representations $(\rho, V_\rho)$ of the parabolic subgroup $P$
\begin{equation}\label{homogen}
\rho \bra {\sc V}_\rho = G \times_P V_\rho~.
\end{equation}

\medbreak
Let $X$ be a compact K\"ahler manifold and $P_G \ra X$  a
holomorphic principal $G$--bundle.
The homogeneous vector bundle ${\sc V}_\rho$ extends to a
holomorphic vector bundle on the associated holomorphic fiber bundle
$M = P_G \times_G F = P_G / P$ by the formula
\begin{equation}\label{ext2}
\wti{\sc V}_\rho \cong P_G \times_P V_\rho \to M = P_G / P~.
\end{equation}
We call $\wti{\sc V}_\rho$ the {\it canonical extension} of
${\sc V}_\rho$~.
With regards to the fundamental group $\Ga = \pi_1(X)$, we suppose
that $M$ has the structure of a generalized flat bundle \cite{KTone}
\begin{equation}\label{flatbundle}
F \hra M = \wti X \times_\Gamma F \oset{\pi} X ~,
\end{equation}
with holonomy $\a : \Gamma \to U$~. Letting
$\o_F$ and $\o_X$ denote the K\"ahler forms of $F$ and $X$ respectively,
the extension to $M$ of the (invariant) K\"ahler form $\o_F$~,
is given by $\wti{\o}_F = p^*\o_F/\a$~, where $p : \wti X \times F \ra F$~,
is the natural projection. Then by \cite{BGKtwo} (Proposition $8.1$),
there exists a family of K\"ahler metrics on $M$ with corresponding weighted
K\"ahler forms
\begin{equation}\label{kform}
\o_{\si} = \pi^* \o_X + \si \wti{\o}_F ~,
\end{equation}
where $\si > 0$ is a constant parameter.

\subsection{The bundle types of the extension on $M$}

Let ${\sc V}_{\rho_i} = U \times_K V_{\rho_i} \to F = U/K$
be homogeneous holomorphic vector bundles with
canonical extensions $\wti{\sc V}_{\rho_i} \ra M$ for $i=1, 2$~.
Further, let
$\sc W_i \ra X$ be holomorphic vector bundles and set
$\sc E_i =  \tpc{\pi^* {\sc W}_i}{\wti{\sc V}_{\rho_i}}$~. We consider
the class of holomorphic vector bundles ${\sc E} \ra M$ given by
proper holomorphic extensions of the form
\begin{equation}\label{ext}
\EBbb ~:~\extn~.
\end{equation}
Such extensions are classified by the $\Ext^1$--functor
(see e.g. \cite{Hart}) which in our case is of the form
\begin{equation}\label{extension5}
\Ext^1_{{\scr O}_M}({\scr E}_2 , {\scr E}_1) \cong H^{0,1}(M,
\schom{{\sc E}_2}{{\sc E}_1} )
\cong H^{0,1} (M, \tpc {\pi^* {\sc W}} {\wti{{\sc V}_{\rho}}} )~,
\end{equation}
where we set ${\sc W} = \schom{{\sc W}_2}{{\sc W}_1}$ and
${\sc V}_\rho = \schom{{\sc V}_{\rho_2}}{{\sc V}_{\rho_1}}$~.
Note that in the latter case we have $\rho = \rho_1 \otimes \rho_2^*$~.

\medbreak
For any holomorphic vector bundle ${\sc W} \to X$ and any homogeneous
vector bundle ${\sc V} \to F$, there is an exact sequence derived from
the Borel--Leray spectral sequence \cite{BGKtwo}
\cite{HZ}~:
\begin{equation}\label{edge}
\begin{aligned}
0 &\ra H^{0,1} (X, \tpc{{\sc W}}{{\sc H}^0 (F, {\sc V}}) )
\osetl{\pi^*} H^{0,1} (M, \tpc{\pi^* {\sc W}}{\wti{\sc V}})
\osetl{\Phi} H^0 (X, \tpc{{\sc W}}{{\sc H}^{0,1} (F, {\sc V}})) \ra \\
&\osetl{d_2} H^{0,2} (X, \tpc{\sc W}{{\sc H}^{0}(F,{\sc V})}) \osetl{\pi^*}
H^{0,2}(M, \tpc{\pi^* {\sc W}}{\wti{\sc V}})~,
\end{aligned}
\end{equation}
where $\pi^*~, ~\Phi$ are the edge homomorphisms.

\medbreak
In the flat case we can say more about the edge map $\Phi$~ \cite{BGKfour}.
Here we use the notation ${\mathbf H}^q (F,{\sc V})$ to indicate
the fact that the fiber cohomologies ${\sc H}^q (F,{\sc V})$ are
flat holomorphic bundles.

\begin{proposition}\label{edge2}~
Suppose that the fiber bundle $F \hra M \to X$ is flat, with holonomy
$\a : \Gamma \to U$~.
For any holomorphic vector bundle ${\sc W} \to X$ and any equivariant
vector bundle ${\sc V} \to F$, we have a short exact sequence

\begin{equation}\label{edge3}
\begin{aligned}
0 &\ra H^{0,1} (X, \tpc{\sc W}{{\mathbf H}^0 (F, {\sc V}}) )
 \osetl{\pi^*}
H^{0,1} (M, \tpc{\pi^* {\sc W}}{\wti{\sc V}}) \\
&\osetl{\Phi}
H^0 (X, \tpc{{\sc W}}{{\mathbf H}^{0,1}(F, {\sc V}})) \ra 0~.
\end{aligned}
\end{equation}
\end{proposition}

This has the following consequence (cf \cite{BGKtwo}~Proposition $7.1$).

\begin{corollary}\label{ExtProp}
Suppose that ${\sc V}$ satisfies the vanishing condition
$$
H^0 (F, {\sc V}) = 0~.
$$

Then the following hold$~:$

\begin{itemize}

\item[(1)]
The holomorphic extensions of the form $\bref{ext}$
are classified by
$$
\Ext^1_{{\scr O}_M} ({\scr E}_2, {\scr E}_1) \cong
H^{0,1}(M, \tpc{\pi^* {\sc W}}{\wti{\sc V}} ) \ocong{\Phi}
H^0 (X, \tpc{\sc W}{{\mathbf H}^{0,1} (F, {\sc V})})~,
$$
and we have
$$
H^0(M, \tpc {\pi^* {\sc W}} {\wti{\sc V}} ) = 0~,
$$
for any holomorphic vector bundle $\sc W$ on $X$~.

\medbreak
\item[(2)]
If $\Gamma = \pi_1 (X)$ acts trivially on
$H^{0,1} (F, {\sc V})$, then
the bundle ${\sc H}^{0,1} (F, {\sc V})$ of fiber co\-ho\-mo\-lo\-gies
is holo\-mor\-phi\-cally trivial and we have the Kunneth formula
$$
H^{0,1}(M, \tpc {\pi^* {\sc W}} {\wti{\sc V}} ) \ocong{\Phi}
\tpc{H^0 (X, {\sc W})}{H^{0,1} (F, {\sc V})}~.
$$
\end{itemize}
\end{corollary}

\subsection{The Kobayashi form of the extension}

In order to describe the representative of the extension class one
needs to construct a right inverse to the edge homomorphism $\Phi$
in Proposition \r{edge2}.
This is done in \cite{BGKfour} and we summarize the necessary results
in the proposition below. It will be useful to keep in mind the
following diagram of holomorphic maps
\begin{equation}\label{maps}
\begin{CD}
\wti X \times F      @>{\ti \pi}>>          \wti X                 \\
      @VVqV                                 @VV{q_0}V              \\
M = \wti X \times_{\Gamma} F    @>{\pi}>>   X                    \\
\end{CD}
\end{equation}
and the relevant cohomology groups as determined by the diagram
\begin{equation}\label{automorphic}
\begin{CD}
H^{0,1} (\wti X \times F, \tpc{\ti q^* {\sc W}}{p^* {\sc V}_\rho})^{\Gamma}
@>{\wti\Phi}>>
\homg{H^{0,1}(F, {\sc V}_\rho)^*}{H^0 (\wti X, q_0^*{\sc W})}{\Gamma} \\
@A{\cong}A{q^*}A        @A{\cong}A{q_0^*}A                       \\
H^{0,1}(M, \tpc{\pi^* {\sc W}}{\wti{\sc V}_\rho})    @>{\Phi}>>
H^0 (X, \tpc{\sc W}{{\mathbf H}^{0,1}(F, {\sc V}_\rho)})~.       \\
\end{CD}
\end{equation}

\begin{proposition}\label{kobform}
With regards to the edge homomorphism $\Phi$ in \bref{edge3}, we
have the following~:

\begin{itemize}

\item[$(1)$]
For a given
$\be_0 \in H^0 (X, \tpc{\sc W}{{\mathbf H}^{0,1}(F, {\sc V}_\rho) })$
there exists a canonical class
$[\bar\be] \in H^{0,1}(M, \tpc{\pi^* {\sc W}}{\wti{\sc V}_\rho})$
such that
$\Phi ([\bar\be]) = \be_0$~.

\medbreak
\item[$(2)$]
The Kobayashi form
$\be \in H^{0,1}(M, \schom{{\sc E}_2}{{\sc E}_1}) \cong
H^{0,1} (M, \tpc {\pi^* {\sc W}} {\wti{\sc V}_\rho})$
of the holomorphic extension \bref{ext}
$$
\EBbb ~:~ 0 \ra \tpc{\pi^* {\sc W}_1}{\wti{\sc V}_{\rho_1}} \lra
{\sc E} \lra \tpc{\pi^* {\sc W}_2}{\wti{\sc V}_{\rho_2}} \ra 0~,
$$
decomposes as
$$
[\be] = \pi^* [\be_X] + [\bar \be]~,
$$
where
$[\be_X] \in H^{0,1} (X, \tpc{{\sc W}}{{\mathbf H}^0 (F, {\sc V}_\rho)})$ and
$\bar \be$ is the right inverse of $\be_0$~.
Thus we have $\Phi ([\be]) = \Phi ([\bar \be]) = \be_0$~.
\end{itemize}
\end{proposition}

Our discussion of extension classes on $M$ leads naturally to the following
definitions of holomorphic objects on $X$~:

\begin{itemize}

\medbreak
\item[$\bullet$]
A {\it holomorphic quadruple}~
$Q = ({\sc W}_1, {\sc W}_2, [\be_X], \be_0)$
is given by two holomorphic vector bundles ${\sc W}_i \to X$,
together with cohomology classes
$$
[\be_X] \in H^{0,1} (X, \tpc{{\sc W}}{{\mathbf H}^0 (F, {\sc V}_\rho )})
~\qquad~,~\qquad~
\be_0 \in H^0 (X, \tpc{{\sc W}}{{\mathbf H}^{0,1} (F, {\sc V}_\rho )})~.
$$
A holomorphic quadruple of the form $Q = ({\sc W}_1, {\sc W}_2, 0, 0)$,
that is $[\be_X] = 0~, ~\be_0 = 0$ is called {\it degenerate}
(see \cite{BGKfour}).

\medbreak
\item[$\bullet$]
A {\it twisted holomorphic triple}~ $T_0 = ({\sc W}_1, {\sc W}_2, \be_0)$
is given by two holomorphic vector bundles ${\sc W}_i \to X$,
together with a holomorphic homomorphism
$$
\be_0 ~:~ {\mathbf H}^{0,1}(F, {\sc V}_\rho)^* \lra
{\sc W} = \schom{{\sc W}_2}{{\sc W}_1}~.
$$
If a basis $\{ {\hat \eta}_j \}$ of $H^{0,1}(F,{\sc V}_\rho)$
is specified, we denote by $T_0$ also the $k$--triple
$T_0 = ({\sc W}_1, {\sc W}_2, \ti \phi)$, where
$\ti \phi = ({\ti \phi}_j)_{j = 1, \ldots, k}$
are the coefficients in the expansion of $\be_0$
(holomorphic triples are considered in \cite{BGPone} \cite{BGKthree}).

\medbreak
\item[$\bullet$]
A {\it twisted $1$--cohomology triple}~
$T_1 = ({\sc W}_1, {\sc W}_2, [\be_X])$
is given by two holomorphic vector bundles
${\sc W}_i \to X$, together with a cohomology class
$$
[\be_X] \in H^{0,1} (X, \tpc{\sc W}{{\mathbf H}^0 (F, {\sc V}_\rho)})~,
$$
classifying a holomorphic extension on $X$ of the form
\begin{equation}\label{extension6}
\WBbb ~:~ 0 \ra \tpc{{\sc W}_1} {{\mathbf H}^0 (F, {\sc V}_\rho)} \lra
\wti{\sc W} \lra {\sc W}_2 \ra 0~
\end{equation}
($1$--cohomology triples are considered in \cite{BGPtwo} ~\cite{DUW} ).


\end{itemize}

Each of the above classes plays a significant role in the dimensional
reduction theory \cite{BGKfour}.
Here we will restrict attention mainly to holomorphic triples and proceed
to state a result which makes use of Corollary \r{ExtProp} and
provides the explicit form of the extension class $\beta$~.

\begin{lemma}\label{BetaLemma2}\cite{BGKfour}~
Suppose that
the homogeneous vector bundle ${\sc V}_\rho$
satisfies the vanishing condition in Corollary $\r{ExtProp}$.

\begin{itemize}

\item[(1)]
Relative to a basis ${\hat \eta}_j = [\eta_j]$ of
$H^{0,1} (F, {\sc V}_\rho)$, the holomorphic triples
$T_0 = ({\sc W}_1, {\sc W}_2, \be_0)$ are of the form
$$
q_0^* ~\be_0 = \summa{j=1}{k}{\ti \phi_j \ot \hat \eta_j}~,
$$
where $\ti \phi =
(\ti \phi_j)_{j = 1, \ldots, k} \in H^0 (\wti X, ~q_0^*{\sc W})^k$
is a $k$--tuple of holo\-mor\-phic sec\-tions.

\medbreak
\item[(2)]
There is a one--one--corres\-pon\-dence between
$k$--tuples $\ti \phi = (\ti \phi_j)_{j = 1, \ldots, k}$
of holo\-mor\-phic sec\-tions and extension classes
$$
[\be] \in
\Ext^1_{{\scr O}_M}({\scr E}_2 , {\scr E}_1) \cong
H^{0,1}(M, \tpc {\pi^* {\sc W}} {\wti{\sc V}} ) \cong
H^0 (X, \tpc{\sc W}{{\mathbf H}^{0,1} (F, {\sc V}_\rho)})~,
$$
given by
$$
q^* \be = \summa{j=1}{k}{\tp{\ti \pi^* \ti \phi_j}{p^* {\eta}_j}}~.
$$

\end{itemize}
\end{lemma}


\section{The perturbation terms associated to a holomorphic extension}

So far we have described how extensions $\EBbb$ in \bref{ext} are
classified by $[\be] \in H^{0,1}(M, \tpc{\pi^* {\sc W}}{\wti{\sc V}_\rho})$
and thanks to Lemma \r{BetaLemma2} we have an explicit form of $\beta$
which will be instrumental in the reduction procedure.
Owing to the generality of our construction, certain technical features
which did not arise in \cite{BGKone} \cite{BGKtwo} now become apparent
and lead to the formulation of the PHE equation.

\medbreak
Henceforth we assume some familiarity with the
differential geometry of operators on K\"ahler manifolds
(references are \cite{Kobthree} \cite{Weil}). In particular, $\lambdas$
will denote the operator of contraction with respect to the K\"ahler form
$\omega_{\si}$ in \bref{kform}.

\subsection{Integration over the fiber}

We define integration over the fiber in the flat fiber bundle $M \to X$~,
\begin{equation}
\pi_* = \fibint ~:~ A^0 (M, ~\tpc{\scend{\pi^* {\sc W}}}
{\scend{\wti{\sc V}_\rho}} )
\lra \emc{X, {\sc W}}~,
\end{equation}
at the level of $\Gamma$--invariant sections
$$
\pi_* = \fibint ~:~
A^0 (\wti X \times F,~ \tpc{\scend{\ti q^* {\sc W}}}
{\scend{p^* {\sc V}_\rho}} )^{\Gamma} \lra
\emc{\wti X, ~q_0^* {\sc W}}^{\Gamma}~
$$
by the formula
\begin{equation}\label{fibint1}
\fibint ~\tp{\ti \pi^* \vphi}{p^* \psi} =
\ov{1}{\Vol (F)}
~\ti \vphi ~\fibint (\psi)  ~\dvol_F =
\frac{1}{\ell! ~\Vol (F)}
~\ti \vphi ~\fibint (\psi)  ~\o_F^\ell~.
\end{equation}
This is well--defined, since the volume form
$\dvol_F = \ov{\o_F}{\ell!}$ is $U$--invariant. 
Here $\Tr$ is induced by the normalized trace on the fiber, that is 
the trace on the bundle $\scend{{\sc V}_\rho}$~. 
We remark that the flatness of the fiber bundle is not necessary
in order to define integration over the fiber.

\medbreak
Observe that the `basic' terms $\be_X$ and $\be_0$ both involve data
on the fiber, holomorphic homomorphisms in the case of $\be_X$ and
holomorphic extensions in the case of $\be_0$~.
There are particular curvature terms
$\lambdas (\be \wedge \be^*)$ and $\lambdas (\be^* \wedge \be)$
which depend on hermitian metrics $h_i$ on ${\sc W}_i$ and the fixed
invariant hermitian metrics $k_i$ on the homogeneous bundles
${\sc V}_{\rho_i}$~.

\begin{lemma}\label{betadecomp}~

\begin{itemize}

\item[(1)]
The endomorphisms~
$-\iota \fibint ~\lambdas (\be \wedge \be^*)
\in \emc{{\sc W}_1}$ ~and~
$\iota \fibint ~\lambdas (\be^* \wedge \be) \in \emc{{\sc W}_2}
$
are non--negative hermitian endomorphisms of ${\sc W}_i~.$

\medbreak
\item[(2)]
If $\fibint ~\lambdas (\be \wedge \be^*) = 0$ or
$\fibint ~\lambdas (\be^* \wedge \be) = 0$~, then $\be = 0$~.

\medbreak
\item[(3)]
For $\be = \pi^* \be_X + \bar\be$ as in Proposition $\r{kobform}$,
we have
$$
\begin{aligned}
\lambdas (\be \wedge \be^*)
&= \lambdas (\pi^* \be_X \wedge \pi^* \be_X^*) +
\lambdas (\bar\be \wedge \bar\be^*)~,     \\
\lambdas (\be^* \wedge \be)
&= \lambdas (\pi^* \be_X^* \wedge \pi^* \be_X) +
\lambdas (\bar\be^* \wedge \bar\be)~.
\end{aligned}
$$
\end{itemize}
\end{lemma}

\medbreak
\begin {definition}\label{vortobstr}
The {\it perturbation terms}  ${\fr d}_i (\be, \si)$
associated to $\be$ are defined by~:
$$
\begin{aligned}
{\fr d}_1 (\be, \si) &= \lambdas (\be \wedge \be^*) -
\tp{\pi^* ~\fibint ~\lambdas (\be \wedge \be^*)}{\wti{\mathbf I}_{1}}~, \\
{\fr d}_2 (\be, \si) &= \lambdas (\be^* \wedge \be) -
\tp{\pi^* ~\fibint ~\lambdas (\be^* \wedge \be)}{\wti{\mathbf I}_{2}}~.
\end{aligned}
$$
\end{definition}

\noindent
The following properties are derived directly from the definition~:

\begin{itemize}

\item[(1)]
$\fibint ~{\fr d}_i (\be, \si) = 0~,$ that is the perturbation terms
vanish under integration over the fiber.

\medbreak
\item[(2)]
For $\be = \pi^* \be_X + \bar\be$ as above, we have
${\fr d}_i (\be, \si) = {\fr d}_i (\pi^* \be_X, \si) +
{\fr d}_i (\bar\be, \si)~,~ i = 1, 2~.$
\end{itemize}

\subsection{The linear maps $\lambda_i$}

Relative to an orthonormal basis $\{ \hat \eta_j \}$ of
$H^{0,1}(F, {\sc V}_\rho)$, we define linear homomorphisms
$$
\lambda_i~:~ \emc{H^{0,1}(F, {\sc V}_\rho)} \to \emc{ {\sc V}_{\rho_i} }~,
$$
by the formulas
\begin{equation}\label{lambdaij}
\lambda_1 (\eta_{ij}) = \frac{1}{\iota} ~\lambdax{F} (\eta_i \wedge \eta_j^*)~,
~\qquad~
\lambda_2 (\eta_{ij}) = \iota ~\lambdax{F} (\eta_i^* \wedge \eta_j)~,
\end{equation}
where $\eta_{ij}$ is the standard basis of $\emc{H^{0,1}(F, {\sc
V})} \cong \gl (k, \CBbb)$~. Since $\eta_{ij}^* = \eta_{ji}$ and
$\lambdax{F} (\eta_i \wedge \eta_j^*)^* = - ~\lambdax{F} (\eta_j
\wedge \eta_i^*)$~, we have $\lambda_i (\xi)^* = \lambda_i
(\xi^*)$~.

\medbreak
There are induced maps on sections
$$
\begin{aligned}
{\lambda_i}_*~&:~A^0 (X, ~\tpc{\scend{{\sc W}_i}}
{\mathbf{End}_{\CBbb}(H^{0,1}(F, {\sc V}_\rho) )} ) \lra A^0 (X,
~\tpc{\scend{{\sc W}_i}} {\mathbf{End}_{\CBbb}({\sc V}_{\rho_i})} )~, \\
\pi^*~&:~A^0 (X, ~\tpc{\scend{{\sc W}_i}}
{\mathbf{End}_{\CBbb}({\sc V}_{\rho_i})}) \lra A^0 (M,
~\tpc{\scend{\pi^* {\sc W}_i}} {\scend{\wti{\sc V}_{\rho_i}}} )~.
\end{aligned}
$$

\medbreak
The main advantage of the maps $\lambda_i$ consists in the fact that
they allow us to express forms like $\lambdas (\bar\be \wedge \bar\be^*)$
in terms of basic data.

\begin{lemma}\label{reduction}~
The forms
$\lambdas (\bar\be \wedge \bar\be^*)$ and
$\lambdas (\bar\be^* \wedge \bar\be)$ are determined by the formulas
\begin{equation}\label{reduction2}
\begin{aligned}
\lambdas (\bar\be \wedge \bar\be^*) &=
\frac{\iota}{\si} ~\pi^* ~{\lambda_1}_* (\be_0 \wedge \be_0^*)~,   \\
\lambdas (\bar\be^* \wedge \bar\be) &=
- ~\frac{\iota}{\si} ~\pi^* ~{\lambda_2}_* (\be_0^* \wedge \be_0)~.
\end{aligned}
\end{equation}
\end{lemma}

\medbreak
For $\be = \pi^* \be_X + \bar\be~, ~\be_0 = \Phi (\bar\be)$
the perturbation terms ${\fr d}_i (\be, \si)$
are now seen to be of the form
$$
{\fr d}_i (\be, \si) =
\pi^* {\fr d}_i (\be_X) + \frac{1}{\si} ~\pi^* {\fr d}_i (\be_0)~.
$$
Observe in particular that the scaling parameter $\si$
for the fiber metric appears as a parameter in the
perturbation terms and this is the reason for the terminology.



\section{The perturbed Hermitian--Einstein equation}

Recalling the holomorphic vector bundle $\sc E \ra M$ in
\bref{ext}, our next step is to describe an integrable unitary
(metric) connection on $\sc E$ and compute its curvature. Then we
will combine this with the background material so far established
and proceed to our intended system of equations.

\subsection{The connection on $\sc E \ra M$ and its curvature}

The holomorphic vector bundle $\sc E$ admits a smooth decomposition
${E} = {E}_1 \oplus {E}_2$~. Relative to this decomposition,
we denote by $\mathbf h$ an invariant hermitian metric on
$\sc E$ of the form
\begin{equation}\label{special}
\mathbf h = {\mathbf h}_1 \oplus {\mathbf h}_2~,
\end{equation}
where ${\mathbf h}_i = h_i' \ot \ti k_i$
on ${\sc E}_i = \tpc{\pi^* {\sc W}_i}{\wti{\sc V}_{\rho_i}}$ is given
by invariant $(basic)$ hermitian metrics $h_i'$ on
$\pi^* {\sc W}_i$ and the extension $\ti k_i$ of $U$--invariant
Her\-mi\-tian--Ein\-stein metrics $k_i$ on ${\sc V}_{\rho_i}~.$

\medbreak
Relative to the smooth decomposition of $\sc E$ and the hermitian metric
$\mathbf h$~, the unitary integrable connection $\mathbf A$ on
$({\sc E}, \mathbf h)$ is given by
\begin{equation}\label{connext1}
\mathbf A = \bmatrix ~\mathbf A_1              & \be
\\ -\be^*               & \mathbf A_2~
\endbmatrix ~,
\end{equation}
where ${\mathbf A}_i$ are the Chern connections of
$({\sc E}_i, \mathbf h_i)~,$ and
$\be \in A^{0,1}(M, \schom {{\sc E}_2} {{\sc E}_1})$
is the representative of the extension class in
$\Ext^1_{{\scr O}_M}({\scr E}_2 , {\scr E}_1)$ relative to
\bref{ext}. A routine calculation (cf e.g. \cite{Kobthree}) shows
that the curvature of $\mathbf A$ has the form
\begin{equation}\label{curvext1}
F_{\mathbf h} = F_{\mathbf A} =
\begin{bmatrix} ~F_{\mathbf h_1} - \be
\wedge \be^*  & D' \be         \\ - D'' \be^*   &   F_{\mathbf
h_2} - \be^* \wedge \be~
\end{bmatrix}~,
\end{equation}
where
\begin{equation}\label{dconnect}
D: A^1 (M, \schom{{\sc E}_1}{{\sc E}_2}) \ra A^2 (M, \schom{{\sc
E}_1}{{\sc E}_2})~,
\end{equation}
is constructed from $\mathbf{A}_1$ and $\mathbf{A}_2$ in the
standard way. Further, we let $D'$ and $D''$ denote the $(1,0)$
and $(0,1)$ components of $D$ respectively, so that $D = D' \oplus
D''$~.

\medbreak
Now for the integrable unitary connection $A_i$
on $({\sc W}_i, h_i)$~, and the Hermitian--Einstein metric
connection ${\wti A}_i$ on $(\wti {\sc V}_{\rho_i}, \ti k_i)$, we have
\begin{equation}
\mathbf A_i = \pi^*A_i \ot \wti{\mathbf I}_{i} +
{\mathbf I}_{i} \ot {\wti A}_i~,
\end{equation}
and the corresponding curvature form of type $(1,1)$ can be expressed as
\begin{equation}\label{curvext2}
F_{\mathbf h_{i}} = F_{\pi^* h_{i}} \ot \wti{\mathbf I}_{i} +
{\mathbf I}_{i} \ot F_{\ti k_i}~,
\end{equation}
where ${\mathbf I}_i = \pi^* \mathbf I_{{\sc W}_i}$ and
$\wti{\mathbf I}_{i} = {\mathbf I}_{\wti{\sc V}_{\rho_i}}$~.

\medbreak
Next, we define the slope of $\sc E$ relative to $\omega_{\si}$
for $\si >0$, by
\begin{equation}\label{muslope}
\lambda = \slos{\EBbb} = \slos{\sc E} = 
\frac{\degs{\sc E}}{\rank (\sc E)}~. 
\end{equation}

\medbreak
We recall that for a compact
K\"ahler manifold $(M, \omega_M)$ and a holomorphic vector bundle
$\sc E \ra M$, the bundle $\sc E$ is {\it stable} ({\it
semistable}) with respect to $\omega_M$, if for any proper
holomorphic subbundle $\sc E' \subset \sc E$ for which $0 < \rank
(\sc E') < \rank(\sc E)$, we have $\mu_{\sc E'} < \mu_{\sc E}$
(respectively, $\mu_{\sc E'} \leq \mu_{\sc E}$). The bundle $\sc
E$ is said to be {\it polystable} if $\sc E$ is a direct sum of
stable bundles of equal slope (see e.g. \cite{Kobthree}).

\medbreak
\begin{definition}\label{PHE}
Relative to  $(M, \o_\si)$ and $\sc E \ra M$ as above, the
{\it perturbed Hermitian--Einstein equation} (the PHE equation)
is defined to be
\begin{equation}\label{defhes}
\iota ~( \lambdas{F_{\mathbf h}} + {\fr d} (\be, \si) ) =
2 \pi ~\lambda ~{\mathbf I}_{{\sc E}}~.
\end{equation}
If ${\fr d} (\be, \si) = 0$, then \bref{defhes}
reduces to the usual Hermitian--Einstein equation \cite{Kobthree}~:
$$
\iota ~\lambdas{F_{\mathbf h}}  =
2 \pi ~\lambda ~{\mathbf I}_{{\sc E}}~.
$$
\end{definition}

\medbreak
\begin{remark}
The term `perturbed' signifies the inclusion of the perturbation
term ${\fr d} (\be, \si)$ in \bref{defhes} which to some extent
accounts for the fact that the usual Hermitian--Einstein equation
does not in general admit a nice decomposition with respect to the
holomorphic fiber bundle $M \ra X$, even in the product case.
The PHE equation necessitates working with a corresponding
restricted type of (poly)stability.
\end{remark}


\section{A dimensional reduction theorem}

\subsection{The calibration condition on the fiber $F$}

We now fix the fiber data as follows. Recall that
on $F$ we are given homogeneous holomorphic bundles
${\sc V}_{\rho_i} \ra F$ associated to complex
representations $(\rho_i~, V_{\rho_i}) \in R(K)$ as in
\bref{homogen}. The degree $\degx{{\sc V}_{\rho_i}}{F}$ and hence the
slope $\mu_{\rho_i} = \sclamj{V}{\rho_i}$ may be computed in terms
of the weights of the representations
$(\rho_i~, V_{\rho_i}) \in R(K)$ by the methods of \cite{BH}.
Further, the representations
$(\rho_i~, V_{\rho_i}) \in R(K)$ are assumed to be sums of
irreducible representations of equal slope. By \cite{Kobthree} IV,
Prop. 6.1, the $K$--invariant hermitian metrics on the irreducible
components of the representation spaces $V_{\rho_i}$ are unique up
to homothety and determine $U$--in\-variant
Her\-mi\-tian--Ein\-stein structures on the corresponding
homogeneous bundles. It follows that the homogeneous bundles
${\sc V}_{\rho_i} \ra F$ have $U$--in\-variant Her\-mi\-tian--Ein\-stein
structures and are therefore polystable. If the representations
$(\rho_i~, V_{\rho_i})$ are irreducible, then the bundles
${\sc V}_{\rho_i}$ are stable (simple) by
\cite{Kobthree} IV, Prop. 6.4 (cf also \cite{Kobtwo} \cite{Ramanan}).
Conversely, if the bundles ${\sc V}_{\rho_i}$ are stable,
the representations $(\rho_i~, V_{\rho_i})$ must be irreducible.
In this case, two homogeneous
bundles ${\sc V}_{\rho_i}$ are isomorphic as holomorphic vector
bundles, if and only if the representations $(\rho_i, V_{\rho_i})$
are equivalent \cite{Ramanan}.

\medbreak
As a consequence, the canonical extensions
$\wti{\sc V}_{\rho_i} \to M$ of the homogeneous
bundles ${\sc V}_{\rho_i} \to F$ admit Hermitian--Einstein structures
\begin{equation}\label{hefiber1}
\iota ~\lambdas F_{\ti k_i} = {2 \pi} ~\ti\mu_{\rho_i}~
\wti{\mathbf I}_{i} ~,
\end{equation}
with constant given by
\begin{equation}\label{hefiber2}
\ti\mu_{\rho_i} = \mu_{\wti{\sc V}_{i}} =
{\frac{\mu_{\rho_i}}{\si}}~.
\end{equation}
It follows that if ${\sc V}_{\rho_i} \to F$ is stable $(\text{simple})$,
then $\wti{\sc V}_{\rho_i} \to M$ is stable $(\text{simple})$.

\medbreak
It is reasonable to require certain {\it calibration conditions}
for the homogeneous vector bundles ${\sc V}_{\rho_i}$~.
In view of the Bott--Borel--Weil theorem \cite{Botttwo}, the representation
theory suggests several possibilities. Here we will assume a slope
condition of the form~
\begin{equation}\label{muslope1}
\mu_{\rho} = \mu_{\rho_1} - \mu_{\rho_2} < 0~.
\end{equation}
It follows that
$H^0 (F, {\sc V}_\rho) = H^0 (F, \schom{{\sc V}_{\rho_2}}{{\sc V}_{\rho_1}}) = 0$
and $V_\rho^K = \homg{V_{\rho_2}}{V_{\rho_1}}{K} = 0$~.
Therefore we are in the situation of Corollary \r{ExtProp}.
Since the cohomology of the irreducible components must now occur in positive
degrees, the Bott--Borel--Weil theorem implies that the corresponding dominant
weights must have Bott index $\geq 1$~.

\begin{remark}
In fact, there are other possible calibration conditions which
could be assumed for other choices of the Bott index and
consequently a different system of equations on $X$ can be
realized \cite{BGKfour}.
\end{remark}

\subsection{The reduction of the PHE equation to the twisted
coupled vortex equations}

\begin{theorem}\label{Main1}~
Let $F \hra M = \wti X \times_\Gamma F \oset{\pi} X$
be a flat holomorphic fiber bundle of compact K\"ahler manifolds.

\medbreak
Suppose that the homogeneous holomorphic bundles ${\sc V}_{\rho_i}$
on $F$ satisfy the calibration condition $\bref{muslope1}~,$
that is $\lam{\rho} = \lam{\rho_1} - \lam{\rho_2} < 0~,$
and therefore $H^0 (F, {\sc V}_\rho) = 0~.$

\medbreak
Consider the proper holomorphic extension
\begin{equation}\label{extension4a}
\EBbb ~:~ 0 \ra \tpc{\pi^* {\sc W}_1}{\wti{\sc V}_{\rho_1}}
\osetl{i} {\sc E} \osetl{p}
\tpc{\pi^* {\sc W}_2}{\wti{\sc V}_{\rho_2}} \ra 0~
\end{equation}
where $\EBbb$ corresponds to the holomorphic triple
$T_0 = ({\sc W}_1, {\sc W}_2, \be_0)~.$

\medbreak
For $\si >0$, let
$$
\lambda = \slos{\EBbb} = \slos{\sc E} = 
\frac{\degs{\sc E}}{\rank (\sc E)}~, 
$$
and define the vortex parameters $\tau_i$ by
$$
\tau_i = \tau_i (\si) = \slos{\sc E} - {\frac{\mu_{\rho_i}}{\si}}~.
$$

\medbreak
Then the following statements are equivalent$~:$

\begin{itemize}

\medbreak
\item[(1)]
There exist invariant hermitian metrics of the form $\mathbf h$ on the
extension bundle $\sc E$ which satisfy the perturbed
Her\-mi\-tian--Ein\-stein equation
\begin{equation}\label{defhes1}
\iota ~( \lambdas{F_{\mathbf h}} + \frac{1}{\si} ~\pi^* {\fr d} (\be_0) ) =
2 \pi ~\lambda ~{\mathbf I}_{{\sc E}}~,
\end{equation}
relative to $(M, \o_\si)$~.

\medbreak
\item[(2)]
There exist hermitian metrics $h_i$ on ${\sc W}_i$ which satisfy
the twisted coupled vortex equations
\begin{equation}\label{genvortex1}
\begin{aligned}
\iota ~\lambdax{X}{F_{h_1}} +
~\frac{1}{\si} ~\fibint ~{\lambda_1}_* (\be_0 \wedge \be_0^*)
&= 2 \pi ~\tau_1 ~\mathbf I_{{\sc W}_1}~,   \\
\iota ~\lambdax{X}{F_{h_2}} -
~\frac{1}{\si} ~\fibint ~{\lambda_2}_* (\be_0^* \wedge \be_0)
&= 2 \pi ~\tau_2 ~\mathbf I_{{\sc W}_2}~.
\end{aligned}
\end{equation}

\medbreak
\item[(3)]
There exist hermitian metrics $h_i$ on ${\sc W}_i$ which satisfy
the twisted coupled multivortex equations
\begin{equation}\label{vortex1}
\begin{aligned}
\iota ~\lambdax{X}{F_{h_1}} + \frac{1}{\si} ~{\Phi_1}
&= 2 \pi ~\tau_1 ~\mathbf I_{{\sc W}_1}~,                        \\
\iota ~\lambdax{X}{F_{h_2}} - \frac{1}{\si} ~{\Phi_2}
&= 2 \pi ~\tau_2 ~\mathbf I_{{\sc W}_2}~,
\end{aligned}
\end{equation}
where $\Phi_i \in \emc{{\sc W}_i}$ are non--ne\-ga\-tive
her\-mi\-tian endo\-mor\-phisms satis\-fying
$$
q_0^* \Phi_1 = \summa{j=1}{k}{\ti \phi_j \circ \ti \phi_j^*}
~\qquad~ ~,~ ~\qquad~
q_0^* \Phi_2 = \summa{j=1}{k}{\ti \phi_j^* \circ \ti \phi_j}~.
$$
Here $k = \dimc H^{0,1} (F, {\sc V}_\rho)$,
and the adjoints $\phi_j^*$ are taken with respect to the
metrics $h_i$~.

\end{itemize}

\medbreak
There is a one--to--one correspondence between solutions in $(1)$
and $(2)$, $(3)$ given by the assignment $h_i \bra h_i' = \pi^*
h_i$~.

\end{theorem}

\medbreak
\begin{remark}
The data in the equations \bref{genvortex1} and \bref{vortex1}
depend only on the associated holomorphic triple
$T_0 = ({\sc W}_1, {\sc W}_2, \be_0)$~.
\end{remark}

\begin{corollary}\label{split}~
If the extension $\EBbb$ is holomorphically split, that is $[\be] = 0~,$
the solutions of the PHE on ${\sc E}$
relative to $(M, \o_{\si})~,$ respectively the corresponding solutions
$(h_1, h_2)$ of the twisted coupled vortex equations $\bref{vortex1}$,
degenerate to solutions of the uncoupled Hermitian--Einstein
equations on each ${\sc W}_i~.$
\end{corollary}

\subsection{Outline of the proof}

The proof of the theorem follows from some technical lemmas which
reflect in part upon the flat structure of \bref{flatbundle}.
We will outline several of the steps involved following \cite{BGKfour}
extending the special cases of \cite{BGKtwo} \cite{GPrfour}.

\medbreak
First of all, we have ~:

\begin{itemize}

\item[(1)]
$\lambdas F_{\pi^* h_i}= \pi^* \lambdax{X} F_{h_i}~;$

\medbreak
\item[(2)]
$\lambdas F_{\ti k_i}= \frac{1}{\si} ~
\wti{\lambdax{F} F_{k_i}}~.$

\end{itemize}

\medbreak
Next, using \bref{curvext1}, \bref{curvext2} and \bref{hefiber1},
the PHE equation is equivalent to the equation
\begin{eqnarray}
&\begin{bmatrix}
~\tp {\pi^* ( \iota ~\lambdax{X} F_{h_1} +
{2 \pi} ~( \ti\mu_{\rho_1} - \lambda ) ~{\mathbf I}_{{\sc W}_1} )}
{\wti{\mathbf I}_{\rho_1}}    & \iota ~\lambdas D' \be     \\
- \iota ~\lambdas D'' \be^*   &
\tp{\pi^* ( \iota ~\lambdax{X} F_{h_2} +
{2 \pi} ~( \ti\mu_{\rho_2} - \lambda ) ~{\mathbf I}_{{\sc W}_2} )}
{\wti{\mathbf I}_{\rho_2}}~
\end{bmatrix}                          \label{curvext3}       \\
&= ~\iota ~
\begin{bmatrix}
~\lambdas ( \be \wedge \be^* ) - {\fr d}_1 (\be, \si)  & 0    \\
0 &  \lambdas ( \be^* \wedge \be ) - {\fr d}_2 (\be, \si)~
\end{bmatrix} ~.                       \notag
\end{eqnarray}
By the definition of the perturbation terms, this last expression equals
\begin{equation*}
\iota ~
\begin{bmatrix}
~\tp{\pi^* ~\fibint ~\lambdas (\be \wedge \be^*)}{\wti{\mathbf I}_{\rho_1}} & 0 \\
0 & \tp{\pi^* ~\fibint ~\lambdas (\be^* \wedge \be)}{\wti{\mathbf I}_{\rho_2}}~
\end{bmatrix} ~.
\end{equation*}
Hence we obtain the equivalent system of equations~:
\begin{equation}\label{lambdabeta}
\lambdas D'\be = 0 ~\qquad~ ~,~ ~\qquad~ \lambdas D''\be^* = 0~,
\end{equation}
and
\begin{equation}\label{genbasic2}
\begin{aligned}
\iota ~\lambdax{X} F_{h_1} - {2 \pi} ~\tau_1 (\si)
~{\mathbf I}_{{\sc W}_1} &=
\iota ~\fibint ~\lambdas ( \be \wedge \be^* )~,              \\
\iota ~\lambdax{X} F_{h_2} - {2 \pi} ~\tau_2 (\si)
~{\mathbf I}_{{\sc W}_2} &=
\iota ~\fibint ~\lambdas ( \be^* \wedge \be )~.
\end{aligned}
\end{equation}

\medbreak
The remainder of the proof deals with some analysis
of $\be$ and showing that the off--diagonal terms in
\bref{curvext3} are zero.
It follows from Lemma \r{BetaLemma2} that we may choose the
smooth decomposition of ${\sc E}$~, such that $q^* \be$ is of the form
$$
q^* \be = \summa{j=1}{k}{\tp{\ti \pi^* \ti \phi_j}{p^* \eta_j}}~,
$$
where $\eta_j \in A^{0,1} (F, {\sc V}_\rho)$ are
$\Delta_{\delbar}$--harmonic $(0,1)$--forms re\-pre\-sen\-ting an
or\-tho\-nor\-mal ba\-sis
of $H^{0,1}(F, {\sc V}_\rho)$~.
Combining this with the Hodge formulas
\begin{equation}\label{Lambdaformula}
\begin{aligned}
\lambdas D' \be - D' \lambdas \be &= \iota ~\delbar^* \be~,  \\
\lambdas \delbar \be - \delbar \lambdas \be
&=
- ~\iota ~{D'}^* \be~.
\end{aligned}
\end{equation}
in \cite{Kobthree}~, we obtain the following equivalent properties
for the metrics $h_i$ on ${\sc W}_i$ (cf \cite{BGKone}).

\begin{lemma}\label{lambdalemma4}~
Suppose that $\Delta_{\delbar} \eta_j = 0$, that is the forms
$\eta_j \in A^{0,1} (F, {\sc V}_\rho)$ are harmonic.
Then we have

\begin{itemize}
\item[(1)]
$\lambdas D' \be = 0~;$

\medbreak
\item[(2)]
$\lambdas D'' \be^* = 0~;$

\medbreak
\item[(3)]
$\delbar^* \be = 0~;$

\medbreak
\item[(4)]
$\Delta_{\delbar} \be = 0$, that is the form
$\be \in A^{0,1} (M, \tpc{\pi^* \sc W}{\wti{\sc V}_\rho})$ is harmonic.

\end{itemize}
\end{lemma}

\medbreak
From the calibration condition \bref{muslope1}, Corollary \r{ExtProp}
and Lemma \r{reduction2} we deduce that
\begin{equation}\label{fibintbeta1}
\begin{aligned}
\frac{1}{\iota} ~\fibint ~\lambdas (\be \wedge \be^*) &=
\frac{1}{\si} ~\fibint ~{\lambda_1}_* (\be_0 \wedge \be_0^*) =
\frac{1}{\si} ~\Phi_1~,  \\
\iota ~\fibint ~\lambdas (\be^* \wedge \be) &=
\frac{1}{\si} ~\fibint ~{\lambda_2}_* (\be_0^* \wedge \be_0) =
\frac{1}{\si} ~\Phi_2~.
\end{aligned}
\end{equation}
Using again the above expression from Lemma \r{BetaLemma2},
the non--negative hermitian endomorphisms $\Phi_i$ of ${\sc W}_i$
admit the expansion
\begin{equation}\label{fibintbeta2}
q_0^* \Phi_1 =
\summa{i,j}{}{\ti \phi_i \circ \ti \phi_j^* ~\langle \eta_i, \eta_j \rangle} =
\summa{j=1}{k}{\ti \phi_j \circ \ti \phi_j^*}
~\qquad~,~\qquad ~
q_0^* \Phi_2 =
\summa{i,j}{}{\ti \phi_i^* \circ \ti \phi_j ~\langle \eta_i^*, \eta_j^* \rangle} =
\summa{j=1}{k}{\ti \phi_j^* \circ \ti \phi_j}
\end{equation}
The Theorem now follows essentially from \bref{genbasic2},
\bref{fibintbeta1} and \bref{fibintbeta2}.

\subsection
{The Reduction Theorem for invariant extensions}

Here we assume the following stronger conditions on the data on the fiber.

\begin{itemize}

\item[(1)]
The representations $(\rho_i, V_{\rho_i}) \in R(K)$ are irreducible.

\medbreak
\item[(2)]
$\mu_{\rho} = \mu_{\rho_1} - \mu_{\rho_2} < 0$~.

\medbreak
\item[(3)]
$H^{0,1} (F, {\sc V}_\rho)^G \neq 0$~.

\end{itemize}

\medbreak
It follows from the Bott--Borel--Weil theorem \cite{Botttwo} and
the multiplicity formulas in \cite{PRV} that the multiplicity of
the trivial representation in $H^{0,1} (F, {\sc V}_\rho)$ is at
most $1$~, that is we have
$H^{0,1} (F, {\sc V}_\rho)^G \cong
H^1 ({\mathfrak p}, {\mathfrak k}_{\CBbb} ; V_\rho) \cong \CBbb$~.
Under the above assumptions, the terms ${\fr d}_i (\be_0)$
vanish and Theorem \r{Main1} takes on a more familiar form.
In fact, we are now essentially in the situation of
\cite[Theorem $8.9$]{BGKtwo}.

\begin{theorem}\label{Main3}
Suppose that the extension $\EBbb$ is invariant, that is
the Kobayashi form $[\be]$ of $\EBbb$ satisfies
$$
{[\beta]} \in \exts{E}{2}{E}{1}{M}_{\mathbf 1} \cong
\tpc{H^0 (X, {\sc W})}{H^{0,1}(F, {\sc V}_\rho)^U}~.
$$

\medbreak
Then the following statements are equivalent:

\begin{itemize}
\item[(1)]
The invariant metric $\mathbf h$ satisfies the Hermitian--Einstein equation
$$
\iota ~\lambdas{F_{\mathbf h}} = 2 \pi ~\lambda ~{\mathbf I}_{{\sc E}}~.
$$

\medbreak
\item[(2)]
The metrics $h_1$ and $h_2$ satisfy
the coupled vortex equations$~:$

\begin{equation}\label{vortex2}
\begin{aligned}
\iota \lambdax{X}{F_{h_1}} + \frac{1}{\sigma} ~\phi \circ \phi^*
&= 2 \pi ~\tau_1 ~\mathbf I_{{\sc W}_1}~,            \\
\iota ~\lambdax{X}{F_{h_2}} - \frac{1}{\sigma} ~\phi^* \circ \phi
&= 2 \pi ~\tau_2 ~\mathbf I_{{\sc W}_2}~.
\end{aligned}
\end{equation}
\end{itemize}
\end{theorem}

\subsection{Examples}

Bott's generalization of the Borel--Weil theorem \cite{Botttwo}
states that for an irreducible $P$--module $(\rho, V_{\rho})$~,
the induced cohomology $H^{0,*} (G/P, {\sc V}_\rho)$ is either
equal to zero or it is an irreducible $G$--module. The theory
underlying the Bott--Borel--Weil (BBW) theorem can be used to
compute examples of irreducible $P$--modules ${\sc V}_{\rho}$
which satisfy the calibration conditions for both of the
reduction theorems as stated above. In principle it seems that a
plentiful supply of such examples can be computed for many types
of the K\"ahler homogeneous space $F= G/P$, in particular for the
case of invariant extensions. The reference \cite{BE} (Chapters
$1$--$5$) outlines a technology for doing this by means of
computational rules. It is based on enumerating the theory of
affine actions of the Weyl group of $\mathfrak g$ and the
Bott--Kostant induction (cf \cite{BGKfour}).

\medbreak
The procedure starts by considering the Dynkin diagram
for a given $\mathfrak g$ where one or more nodes $\bullet$ are
replaced by a crossed node $\times$ when there is a non--parabolic
simple root. In this way the Dynkin diagram for $F = G/P$ is
obtained. For instance, a maximal parabolic $\mathfrak p$
subalgebra (as in the case of the compact irreducible Hermitian
symmetric spaces) admits a single crossed node and a Borel
subalgebra $\mathfrak b$ has crosses through every node as is the
case for the full flag manifold $G/B$ over $\mathbb C^{\ell +1}$~.
The weights of the representations are exhibited by such diagrams
(see below) by placing (integer) coefficients over each node in
accordance with certain rules. If the aim is to obtain a
$1$--dimensional irreducible $\mathfrak g$--module, we would
select suitable node coefficients for the diagram such that on
taking a single affine Weyl group reflection over the appropriate
crossed node leads to zeros over each node in the diagram and
hence this selection corresponds to the trivial $\mathfrak g$--module
provided by the BBW theorem.

\begin{example} Let $F = \mathbb CP^4$ and take $\sc V_\rho$ to be the
irreducible $P$--module $\O_F^1$~. Starting from the corresponding
Dynkin diagram~ $\ \ldyncbbb{-2}{1}{0}{0}$~ and then taking a
single affine reflection, the trivial module is obtained. Thus
$H^{1}(F, {\sc V}_\rho)\cong \mathbb C$~, and the cohomology in all
other degrees is zero by the BBW Theorem.

\medbreak
The dual module ${\sc V}_\rho^* \cong \sc T^{1,0}_F$ has the
corresponding diagram~ $\ \ldyncbbb{1}{0}{0}{1}$~. That for the
canonical line bundle $\sc K_F = \O_F^4$ is ~$\
\ldyncbbb{-5}{0}{0}{0}$~ from which ${\sc V}_\rho \ot {\sc K}_F$ is
represented by ~$\ \ldyncbbb{-7}{1}{0}{0}$~. On taking four affine
reflections on the latter we obtain ~ $\ \ldyncbbb{1}{0}{0}{1}$~.
Thus Serre--duality and the BBW theorem imply that the irreducible
$G$--module
$$
H^0 (F, {\sc V}_\rho^*) \cong H^4 (F, {\sc V}_\rho \ot {\sc K}_F)
\cong {\fr g}~:~~ \ldyncbbb{1}{0}{0}{1}
$$
and the cohomology in all other degrees is zero,
so we have $H^0(F, {\sc V}_\rho^*)^U = H^1 (F, {\sc V}_\rho^*)^U = 0$~.
\end{example}

\medbreak
\begin{example}
Here we take $F$ to be the $9$--dimensional partial flag manifold
over $\CBbb^5$ whose compact representation is
$$
F \cong SU (5)/ S (U (1) \times U (2) 
\times U (1) \times U (1))~.
$$
It is an example of a homogeneous
K\"ahler manifold which is not a symmetric space. Consider the
irreducible $P$--module ${\sc V}_\rho$ as represented by ~
$\ \ldyncbcc{-2}{1}{0}{0}$~. A single affine reflection leads to
$H^{1}(F, {\sc V}_\rho) \cong \mathbb C$ and zero cohomology in all
other degrees. As for the dual module $\sc V_\rho^*$, the diagram is ~
$\ \ldyncbcc{1}{1}{-1}{0}$~ which can be seen to correspond to a
singular weight and hence $H^q(F, {\sc V}_\rho^*) = 0$, for all
$q \geq 0$~.
\end{example}

\medbreak
\begin{remark}
We remark that other types of examples can be formulated following
e.g. \cite{MS} Theorem B.
\end{remark}


\section{The moment map and the PHE equation}

Given a hermitian metric $h$ on $E$~, the space ${\sc C} (E)$ of
integrable $\delbar$--operators, for which $\delbar^2 = 0$~,
corresponds bijectively to the space ${\sc A} (E, h)$ of unitary
integrable connections whose curvature satisfies $F^{0,2}_{h} = 0$~.
Here $E$ denotes the underlying smooth vector bundle of $\sc E$~.
Thus each element $\delbar_E \in \Ce (E)$ defines a
unique holomorphic structure ${\sc E} = (E, \delbar_E)$ on $E$~,
for which it is the canonical $\delbar$--operator.
The complex gauge group $\Aut (E)$ acts on $\Ce (E)$ via the action
$g(\delbar) = g \circ \delbar \circ g^{-1}$~, for $g \in \Aut (E)$~.
For our purpose, we restrict attention to the {\it unitary} gauge
(sub)group denoted by $\G$~.
The quotient $\Ce (E)/{\G}$ is the space of equivalence classes of
integrable holomorphic structures on $E$ up to unitary equivalence.
The Lie algebra of $\G$ is given by
$\Lie (\G) \cong \End_s (E)$, where $\End_s (E)$ is the Lie algebra
of global skew--hermitian endomorphisms of $E$~.
Background references to this section are \cite{AB}
\cite{Donone} \cite{Kobthree}.

\subsection{The restricted gauge group and the moment map}

Since $M$ is K\"ahler, the inner product
\begin{equation}
\langle \a_1,\a_2 \rangle = \frac{1}{\iota (m-1)!\Vol(M)}
\int_M \Tr ~(\a_1 \wedge \a_2^*) \wedge \o_{\sigma}^{m-1}~,
\end{equation}
for $\a_1,\a_2 \in T_{\delbar}~\Ce (E) \cong  A^{0,1}(M,
\End_s (E))$, induces a K\"ahler structure on $\Ce (E)$ where
the K\"ahler form $\o$ is defined by $ \o(\a_1, \a_2) =
\im~\langle \a_1,\a_2 \rangle$~. The standard action of $\G$ on
$\Ce (E)$ preserves $\o$ and induces an associated
equivariant moment map
\begin{equation}\label{moment0}
\begin{aligned}
\nu = \nu (\G)~: ~ \Ce (E) &\lra
\Lie (\G) \subset \Lie (\G)^* \cong L^2({\Lie (\G)})  \\
\delbar &\mapsto \lambdas{F_h} ~,
\end{aligned}
\end{equation}
where $\delbar$ corresponds to the unitary integrable connection
$(A, h)$~. This moment map is determined up to a constant
in the center of the Lie algebra and may also be written as
\begin{equation}\label{einsteinmoment}
\nu (\delbar) = \lambdas{F_{h}} ~+~ 2 \pi \iota ~\lambda
~{\mathbf I}_{\sc E}~.
\end{equation}
Note that $\nu^{-1} (0)$ is empty unless $\lambda = \mu_{E}$~,
the slope of $E$~.

\medbreak
In this section we assume that the representations $(\rho_i,
V_{\rho_i})$ are irreducible. We consider the subspace ${\sc A}
(\EBbb, {\mathbf h}) \subset {\sc A} (E, {\mathbf h})$ of unitary
integrable connections ${\mathbf A}$ of the form $({\mathbf A}_1,
{\mathbf A}_2, \be)$~, as in \bref{connext1}, where $\mathbf h$ is
a (fixed) special metric $\mathbf h$ as in
\bref{special}.
Let $\Ce (\EBbb) \subset \Ce (E)$ be the subspace of holomorphic
structures determined by ${\sc A} (\EBbb, {\mathbf h})$~.
The elements $\delbar_E \in \Ce (\EBbb)$ determine a holomorphic
structure on the extension $\mathbb E$ in \bref{extension4a}, that
is $\delbar_E$ is of the form
\begin{equation}\label{holomoment}
\delbar_E =
\begin{bmatrix}
~\delbar_{\pi^* W_1} \otimes \wti{\mathbf{I}}_1 + \mathbf{I}_1 \ot
\delbar_{\wti{\sc V}_1} & \beta    \\ 0 & \delbar_{\pi^* W_2} \ot
\wti{\mathbf{I}}_2 + \mathbf{I}_2 \ot \delbar_{\wti{\sc V}_2}~
\end{bmatrix}~.
\end{equation}
We further consider the subspaces
$\Ce_0 (\EBbb) \cong {\sc A}_0 (\EBbb, {\mathbf h})$~, consisting of
the elements in $\Ce (\EBbb) \cong {\sc A} (\EBbb, {\mathbf h})$
such that $\be$ is $\Delta_{\delbar}$--harmonic, that is
$\delbar^* \be = 0~.$
Then ${\sc A}_0 (\EBbb, {\mathbf h})$ admits a mapping
\begin{equation}\label{fiber}
{\sc H}_{\delbar}^{0,1} (M, \schom{{\sc E}_2}{{\sc E}_1}) \lra
\A_0 (\EBbb, \mathbf h) \osetl{\Pi} \A (W_1, h_1) \times \A (W_2, h_2)~,
\end{equation}
where the dimension
$\dimc {\sc H}_{\delbar}^{0,1} (M, \schom{{\sc E}_2}{{\sc E}_1})$
is upper--semicontinuous as a function on the base.

\medbreak
We specify a subgroup $\G_0 \subset \G$ which acts symplectically
on $\Ce (E)$ and on $\Ce_0(\EBbb)$ via restriction to the latter.
For $u_i \in \G_{W_i}$, the subgroup $\G_0$ is defined by
\begin{equation}\label{gaugetwo}
\G_0  = \{
\begin{bmatrix} ~
\pi^* u_1 \ot \wti{\mathbf{I}}_{1} & 0 \\ 0 &
\pi^* u_2 \otimes \wti{\mathbf{I}}_{2}~
\end{bmatrix}
\}
\cong  \G_{W_1} \times  \G_{W_2} \subset \G ~.
\end{equation}
The subgroup $\G_0$ leaves $\Ce_0(\EBbb)$ invariant and fixes the
holomorphic structures on the fiber. In fact, $\G_0$ is the
maximal subgroup of $\G$ with this property, since by the
irreducibility of the $V_{\rho_i}$ there are no non--trivial
$U$--equivariant gauge transformations on the homogeneous bundles
${\sc V}_{\rho_i}$~, that is
$\emg{{\sc V}_{\rho_i}}{U} \cong \emc{V_{\rho_i}}^K \cong
{\CBbb} \cdot \Id$~.

\medbreak
On smooth elements (relative to
$\mathbf h = \mathbf h_1 \oplus \mathbf h_2$ as above),
we have a commutative diagram
\begin{equation}\label{orthomoment2}
\CD
\Lie (\G) @>\subset>> \Lie (\G)^* \cong L^2({\Lie (\G)}) \\
 @AA\pi^*A    @VP_0VV      \\
\Lie (\G_0) @>\subset>> \Lie (\G_0)^* \cong L^2({\Lie (\G_0)}) \\
\endCD
\end{equation}
where $P_0$ denotes orthogonal projection and for $a \in \Lie
(\G_0)$ the following relationship is satisfied~:
\begin{equation}
\langle
P_0(\lambdas{F_{\mathbf h}}) ~,~ a \rangle =
\frac{\iota}{n!\Vol(X)} \int_X \Tr~
(P_0(\lambdas F_{\mathbf h}) \circ a^*)~\o_{X}^n~.
\end{equation}

\medbreak
Observing that the projection $P_0$ is essentially given by
integration over the fiber, we obtain the main result concerning
the moment map interpretation
of the PHE equation.

\begin{theorem}\label{mainmoment}\cite{BGKfour}~
With regards to the inclusion
$j : \Ce_0 (\EBbb)  \hookrightarrow \Ce (E)$~, consider the map
$$
\nu_0  : \Ce_0 (\EBbb) \lra \Lie(\G_0)^* ~,
$$
as defined by $\nu_0 = P_0 \circ \nu \circ j$~. Then the
following hold ~:

\begin{itemize}
\item[(1)]
The map $\nu_0$ is a moment map for the action of $\G_0$ on
$\Ce_0 (\EBbb)$~.

\medbreak
\item[(2)]
The following diagram commutes with respect to the inclusion of
smooth elements$~:$
$$
\begin{CD}
\Ce (E)
 @>\nu >>  \Lie (\G) @>\subset >> \Lie (\G)^* \cong L^2({\Lie (\G)}) \\
 @AAjA                @AA\pi^*A    @VP_0VV      \\
\Ce_0 (\EBbb)
 @>\nu_0 >>  \Lie (\G_0) @>\subset>> \Lie (\G_0)^* \cong L^2({\Lie (\G_0)}) \\
\end{CD}
$$

\medbreak
\item[(3)]
The PHE equation
$$
\iota (\lambdas F_{\mathbf h} ~+~ {\fr d} (\be, \si)) =
2 \pi ~\lambda ~\mathbf I_{\sc E}~,
$$
is equivalent to the $\G_0$--moment map equation $\nu_0(\delbar)=0$~.

\end{itemize}
\end{theorem}



\end{document}